\documentclass[article,ij4uq]{ij4uq}      

\usepackage[hang]{footmisc}
\fancypagestyle{plain}{%
  \fancyhf{}
  \fancyhead[R]{\small {\it International Journal for Uncertainty Quantification}, x(x): \thepage--\pageref{LastPage} (\myyear\today)}
  \fancyfoot[R]{\small\bf\thepage }
  \fancyfoot[L]{\fottitle}
  }
\renewcommand{\dmyy}{17}
\renewcommand{\myyear}{2017}
\renewcommand{\today}{}

\usepackage[english]{babel}
\usepackage[utf8]{inputenc}
\usepackage{mathabx}
\usepackage{graphicx,epstopdf}
\usepackage{cleveref}

\usepackage{tikz}
\usepackage{pgfplots}
\pgfplotsset{compat=1.13}
\usepgfplotslibrary{external}
\tikzset{external/system call={lualatex -shell-escape -interaction=batchmode -jobname "\image" "\texsource"; convert -density 600 "\image.pdf" "\image.eps"}}
\tikzset{/pgf/images/external info}
\tikzset{/pgf/images/include external/.code={\includegraphics[width=\pgfexternalwidth,height=\pgfexternalheight]{#1.eps}}}
\tikzsetexternalprefix{./}
\DeclareMathOperator*{\argmin}{arg\,min}

\renewcommand{\vec}[1]{\bm{#1}}
\newcommand{\mat}[1]{\bm{#1}}
\newcommand{\kernel}{k}

\newcommand{\R}{\mathbb{R}}
\newcommand{\N}{\mathbb{N}}
\newcommand{\dom}{\mathcal{D}}

\begin{document}
\volume{Volume x, Issue x, \myyear\today}
\title{Ensemble Kalman Filters for reliability estimation in perfusion inference}
\titlehead{EnKF for reliability estimation in perfusion inference}
\authorhead{P.~Zaspel}

\corrauthor[1]{Peter Zaspel}
\corremail{peter.zaspel@unibas.ch}
\corraddress{University of Basel, 4051 Basel, Switzerland}

\dataO{mm/dd/yyyy}
\dataF{mm/dd/yyyy}

\abstract{We consider the solution of inverse problems in dynamic contrast--enhanced imaging by means of Ensemble Kalman Filters. Our quantity of interest is blood perfusion, i.e.~blood flow rates in tissue. While existing approaches to compute blood perfusion parameters for given time series of radiological measurements mainly rely on deterministic, deconvolution--based methods, we aim at recovering probabilistic solution information for given noisy measurements. To this end, we model radiological image capturing as sequential data assimilation process and solve it by an Ensemble Kalman Filter. Thereby, we recover deterministic results as ensemble--based mean and are able to compute reliability information such as probabilities for the perfusion to be in a given range. Our target application is the inference of blood perfusion parameters in the human brain. A numerical study shows promising results for artificial measurements generated by a Digital Perfusion Phantom.}

\maketitle

\pgfplotsset{every axis legend/.append style={font=\footnotesize}}

\section{Introduction}\label{sec:introduction}
Medical imaging by x-rays, \textit{magnetic resonance imaging} (MRI) and \textit{computed tomography} (CT) has considerably changed medical diagnosis throughout the last decades.
Often, \textit{contrast agents}, i.e.~specific liquid chemicals, are injected into the patients blood circulation during the imaging process. This leads to \textit{contrast--enhanced images} of higher contrast in some regions of the human body. In this work, we study inverse problems for a specific MRI or CT imaging task. That is, we aim at recovering a quantity of interest in medical imaging, which is derived by \textit{dynamic} contrast--enhanced (DCE) imaging. In dynamic contrast--enhanced imaging, a time-dependent series of radiological images of a part of the patients body (e.g.~the brain) is taken immediately after injecting a contrast agent into the patient's blood circulation. By observing the time-dependent concentration evolution of the contrast agent inside the patient's tissue, it is possible to recover information about the blood flow rates, i.e.~the \textit{perfusion}. 

The outcome of the image acquisition process is a time-discrete series of tree-dimensional (space-discrete) concentration images $\vec{c}$ of a part of the patient body. The actual perfusion evaluation is a post-processing step, being preceded by image de-noising and motion compensation. Currently, blood perfusion is computed independently per discrete tissue volume element, i.e.~\textit{voxel}, thus spatial information is mostly neglected. Variants of the \textit{indicator-dilution theory} \cite{Brix2010,Fieselmann2011,Sourbron2014} describe the concentration of a contrast agent in tissue at a given point in time as the result of a convolution in time of the (known) time-dependent arterial or blood circulation inflow concentration $c_{art}$ with an unknown tissue-dependent kernel function $\vec{\kernel}$. \textit{Blood perfusion} is computed as a weighted maximum or point evaluation of the unknown kernel function.

Current state of the art methods aim at recovering the unknown time-dependent kernel function for given discrete measurements of the contrast agent in tissue. The kernel function is either modeled as a parametrized analytic function \cite{Tofts1997,Parker2006,Fieselmann2011} or discretized as a fully unknown function \cite{Brix2010,Ostergaard1996,Fieselmann2011}. Then, most approaches rely on a deterministic reconstruction of the kernel function, involving the solution of a deconvolution problem with regularization.

One drawback of the use of motion compensation, de-noising and deterministic deconvolution lies in the loss of information on the \textit{quality} of a computed solution. That is, probabilistic information about the measurement accuracy and errors in space and time with their influence on the exactness on the computed quantity of interest are neglected or even lost. 
          
In this work, we propose an approach to infer perfusion in the discussed application case while keeping the probabilistic information on the solution. Thereby, we overcome the discussed drawback of knowledge loss. To achieve this, we model the inference problem as a sequential data assimilation problem: First, the unknown kernel function $\vec{\kernel}$ is described as an unknown system state, for which a predictive time-discrete stochastic system state model is introduced. In this model, the kernel function is represented as a random variable. Then, a time-discrete stochastic observation model describes the relationship of the current approximation of $\vec{\kernel}$ and the noisy measurements delivered by medical imaging. Finally, the well-known Ensemble Kalman Filter (EnKF) \cite{Evensen2009,Stuart2010,Iglesias2013,Ernst2015} is used to compute an ensemble--based approximation of the posterior probability density function (PDF) of $\vec{\kernel}$ given the system state model and the (noisy) measurements. Based on this PDF, means, cumulative distribution functions, etc.~can be computed. 
The whole sequential data assimilation methodology is applied to (noisy) artificial measurement data generated by a \textit{Digital Perfusion Phantom} \cite{Riordan2011,Pianykh2012,Manhart}, i.e.~a forward model describing the mapping of perfusion information to medical images. Note that we stick to the use of a Digital Perfusion Phantom, instead of using exemplary patient data since, first, it allows to artificially create arbitrary amounts of radiological images and, second, a scanner- and patient-independent way to analyse perfusion estimation methods is highly desired in radiology. Certainly, applying the proposed methodology (together with radiologists) to real patient data is future work.
   
The data assimilation problem that we will model in Section~\ref{sec:numerical_methods} will stick to Gaussian random fields and a linear forward model. This is a strong simplification, allowing to analytically solve the data assimilation problem by the original \textit{Kalman filter} \cite{Reich2015}, However, we use the EnKF, which is a generalization of this Kalman filter, usually being applied in a non-linear, non-Gaussian setting. For more details on the connection between EnKF and the Kalman filter and a convergence analysis in the linear Gaussian case see \cite{Majda2018,Tong2018}, while EnKFs for inverse problems are discussed in \cite{Iglesias2013,Schillings2017,Schillings2018} and many extensions and alternatives for the EnKF are, e.g., developed in \cite{Anderson2001,Anderson2003,Tippett2003,Zhang2009,Reich2013,Reich2013a}. We decided to use the EnKF here, since we consider this work as a starting point for much more involved approaches for the prediction of perfusion. In fact, the rather simple linear forward model from the indicator dilution theory should be replaced by more complex or even PDE-based models, which will certainly no longer be linear. Moreover, we expect that a more involved evolution model, cf.~Section~\ref{sec:evolutionModel}, with non-Gaussian noise might become valid in real application cases. Therefore, we here already introduce the more involved EnKF framework, while considering other forward models and non-Gaussian noise as future work. 

To the best of our knowledge, we consider the discussed work to be a new contribution to the field. Nevertheless, there has been previous work on the use of Ensemble Kalman Filters in the application scenario. In \cite{Naevdal2016}, the authors concentrate on the introduction of a tissue model that includes space-dependent information. To achieve this, a blood flow model is combined with an EnKF. Preliminary results for this approach are given. In contrast, we focus here directly on the mathematical setting based on the indicator-dilution theory that is well-known and, therefore, well accepted by radiologist. Hence, our methodology is considered as an extension to the existing standard methodology introducing the opportunity to derive statistical information on the computed solution. In addition to the different objective compared to \cite{Naevdal2016}, we also perform a large number of parameter studies and convergence tests, which are crucial to understand the properties of the method.

This article is organized as follows. In \Cref{sec:problem_formulation}, we give a mathematical model for the radiological imaging and perfusion extraction mechanism. \Cref{sec:numerical_methods} outlines our numerical approach based on sequential data assimilation using EnKF. Numerical results are given in \Cref{sec:numerical_results} while \Cref{sec:summary} summarizes the discussed work.

\section{Modeling radiological imaging and perfusion extraction}\label{sec:problem_formulation}
In the following, we start by giving an abstract model for the transport of contrast agent. Then, measurements by e.g.~MRI are abstractly modeled. A concrete model for the contrast agent distribution is given by the indicator-dilution theory. Finally, our quantity of interest, i.e.~blood perfusion, is introduced and the deterministic inference problem is summarized.
\subsection{Abstract model for contrast agent transport}
Shall $\dom_{tiss}\subset\R^3$ be the tissue domain in the human body for which we want to derive information by dynamic contrast--enhanced imaging. 
We study contrast agent transport / concentrations in a time interval $[0, T]\subset\R$ with $T$ being the final time. The inflow concentration of the contrast agent (at some arterial inlet) is a function
$c_{art}:[0,T] \rightarrow \R_{\geq 0}\,.$ 
The time-continuous contrast agent concentration in tissue can be modeled as a function
$c: \dom_{tiss}\times [0,T] \rightarrow \R_{\geq 0}\,.$
Both are related to each other by an (unknown) operator $\mathcal{B}$, with
\begin{equation}\label{eq:space_continuous_transport_model}c(\cdot,t) = \mathcal{B}[c_{art}](t)\end{equation}
that models the function of the human body with respect to contrast agent transport.

\subsection{Measuring contrast agent concentration in tissue}
Appropriate measurement devices (CT, MRI, \ldots) usually have a cuboidal measurement domain. Therefore, we start by limiting $\dom_{tiss}$ to
$\dom_{meas}= \bigtimes_{d=1}^3 [0, a_d]$
with $\vec{a}=(a_1, a_2, a_3)^\top\in\R^3$ describing the size of the measurement domain. For simplicity, we assume the measurement domain and the area of interest to match exactly, i.e.~$\dom_{meas} = \dom_{tiss}$, excluding cases in which some part of the measurement domain does not contain valid tissue. Moreover, $\dom_{meas}$ is simplified as being stationary in time, i.e.~the measurement device (or the patient) does not move or movements are considered as measurement error. 

The finite spatial resolution $\vec{N}_{\dom}\in\N^3$ of the measurement device leads to a decomposition of $\dom_{meas}$ into $N_{voxel} = \prod_{d=1}^3 N_{\dom}^{(d)}$ volume elements or \textit{voxels} of volume $V_{voxel} = \prod_{d=1}^3 a_d / N_\dom^{(d)}$ for which we obtain averaged (constant) measurements. 
We introduce a measurement operator $\Psi$ that 
gives for a given exact contrast agent concentration $c$ and a chosen point in time $t\in [0,T]$ a measurement vector $\vec{c}(t)\in\R_{\geq 0}^{N_{voxel}}$ as
$$ \vec{c}(t) = \Psi[c](t) := \Theta[c](t) + \mathcal{E}[c](t)\,.$$
Here, $\Theta$ is a noise-free measurement-operator and is usually a volumetric average over each voxel being equivalent to a piece-wise constant approximation in space.
$\mathcal{E}$ abstractly models a (potentially non-linear) additive error (noise, movements, technical problems, \ldots).  

To reflect time-discrete measurements, we introduce $N_{obs}$ ordered, pair-wise different discrete observation times $t_i^{obs}\in [0,T]$, $i\in \{1, \ldots, N_{obs}\}$, at which measurements or observations are done, giving the observation matrix $\mat{C}: = (c_{ji})_{j=1,\ldots, N_{voxel}, i=1,\ldots, N_{obs}} $ composed of observation vectors $\vec{c}_i$ as $\mat{C}= (\vec{c}_1| \ldots |\vec{c}_{N_{obs}}) $ with
\begin{equation}\label{eq:measurement_operator_for_observations}\vec{c}_i = \Psi[c](t_i^{obs}) := \Theta[c](t_i^{obs}) +\mathcal{E}[c](t_i^{obs}) \,.\end{equation}

\subsection{Contrast agent transport model following the indicator-dilution theory}
The in\-di\-ca\-tor-dilution theory (IDT) \cite{Brix2010} provides a model for the time evolution of the contrast agent concentration in a reference voxel $\dom_{voxel}\subset\dom_{tiss}$ with volume $V_{voxel}$, given the arterial inflow $c_{art}$. While, in this standard model, the contrast agent's concentrations are assumed to be constant in each voxel, we first want to formulate the IDT as a space-continuous model and then move over to a discrete description as consequence of a measurement process. Our continuous version of the indicator-dilution theory--based transport model replaces $\mathcal{B}$ in \cref{eq:space_continuous_transport_model} with the model operator $\mathcal{B}_{IDT}$ given via
\begin{equation}\label{eq:continuous_idt}c(\vec{x},t) =  \mathcal{B}_{IDT}[c_{art},\kernel](t) := \int_{0}^T c_{art}(\tau) {\kernel}(\vec{x},t-\tau) d\tau\,,\quad (\vec{x},t)\in \dom_{tiss}\times [0,T]\,.
\end{equation}
Kernel ${\kernel}:\dom_{tiss}\times [0,T] \rightarrow\R$ fully characterizes the properties of the tissue at point $\vec{x}$. In order to have an well-defined integrand, we assume $\kernel(\cdot,t)=0$ for $t<0$. Note that the model operator $\mathcal{B}_{IDT}$ is actually independent of the spatial position.  

We now apply the measurement operator $\Psi$ to \cref{eq:continuous_idt} obtaining
\begin{equation*}
\begin{aligned}
\Psi[c](t) &= \Theta\left[\mathcal{B}_{IDT}[c_{art},\kernel]\right](t) + \mathcal{E}\left[\mathcal{B}_{IDT}[c_{art},\kernel]\right](t)\\
 & = \int_0^T c_{art}(\tau)\vec{\kernel}(t-\tau) d\tau + \mathcal{E}\left[\mathcal{B}_{IDT}[c_{art},\kernel]\right](t)\,,
\end{aligned}
\end{equation*}
where $\vec{\kernel} = (\kernel_1,\ldots ,\kernel_{N_{voxel}})^\top$ is a vector of univariate kernel functions $k_j:[0,T]\rightarrow \R$. Since we are interested in time-discrete observations, we limit our discussion to observation times $t_i^{obs}$ yielding
$$\vec{c}_i = \Psi[c](t_i^{obs}) = \int_0^T c_{art}(\tau)\vec{\kernel}(t_i^{obs}-\tau) d\tau + \vec{e}_i(c_{art},\kernel)\,,$$
with the abbreviation $\vec{e}_i(c_{art},\kernel) :=  \mathcal{E}\left[\mathcal{B}_{IDT}[c_{art},\kernel]\right](t_i^{obs})$.
For a single voxel\linebreak $j\in\{1,\ldots, N_{voxel}\}$, we obtain
$$c_{j,i}^{obs} =  \int_0^T c_{art}(\tau)\kernel_j(t_i^{obs}-\tau) d\tau + {e}_{j,i}(c_{art},\kernel)\,.$$
In case of ${e}_{j,i}(c_{art},\kernel) = 0$, this boils down to the classical indicator-dilution-theory model given on a reference voxel $j$. Obviously, this model is independent of the spatial position of the voxel $j$. The classical theory further introduces a mean density $\rho_j\in\R_{\geq 0}$ in a voxel $j$, which becomes of interest in the following subsection.

\subsection{Perfusion}
The inference task discussed in this article is to compute a time-stationary perfusion (blood flow) information $\vec{p}\in\R^{N_{voxel}}$ given the (assumed to be exactly known) inflow concentration $c_{art}$ and the observation matrix $\vec{C}$. Formally, the blood perfusion in a given voxel $j$ can be evaluated as quantity of interest of the computed response function $\kernel_j$ as
$${p}_j := {p}(\kernel_j) := \frac{1}{\rho_{j}} k_j(0)\,.$$
From a mathematical point of view, this quantity has nice properties, since it is just a point evaluation of the response function. In practice \cite{Brix2010}, perfusion is however often evaluated as
\begin{equation*}
\tilde{p}_j := \tilde{p}(\kernel_j) := \frac{1}{\rho_{j}} \max_{t\in[0,T]}\kernel_j(t)\,.
\end{equation*}
For simplicity and since we use just artificial input data, we stick to the first version of this quantity of interest.

\subsection{Deterministic inference problem}\label{sec:full_deterministic_problem}
To summarize this section, we formulate the deterministic problem that we aim to solve: For given measurement time $T\in\R$, arterial inflow $c_{art}:[0,T]\rightarrow \R_{\geq 0}$, measurement/observation times $t_1^{obs}<t_2^{obs}<\ldots < t_{N_{obs}}^{obs}$ and observation matrix $\mat{C}$ or vectors $\vec{c}_i\in\R^{N_{voxel}},\, i\in\{1,\ldots,N_{obs}\}$, we aim at computing a vector $\vec{\kernel}=(\kernel_1,\ldots,\kernel_{N_{voxel}})^T$ of kernel functions ${\kernel}_j:[0,T]\rightarrow\R$ and the derived quantity of interest $\vec{p}=(p_1, \ldots , p_{N_{voxel}})^\top$ with $p_j = \kernel_j(0) / \rho_j$ such that 
\begin{equation}\label{eq:deterministic_observation_equation}
c_{j,i} \approx \int_0^T c_{art}(\tau)\kernel_j(t_i^{obs} - \tau) d\tau + e_{j,i}(c_{art},\kernel)\,,\quad\quad j\in\{1,\ldots,N_{voxel}\},\,\, i\in\{1,\ldots,N_{obs}\}\,.
\end{equation}
Clearly, this problem is underdetermined with the given requirements. Furthermore, we have not specified the nature of the error term, yet. This is why we used the notion ``$\approx$''. A much clearer idea of the concept of a \textit{solution} to this problem is given in the next section, where we reformulate the problem as Bayesian sequential data assimilation problem.
\section{Numerical approach by sequential data assimilation}\label{sec:numerical_methods}
In this section, we first introduce a discretization for the model discussed in the last section. This is necessary, since we will use its discretized version in context of sequential data assimilation, afterwards. An approximation to the solution of the assimilation problem is derived by the Ensemble Kalman Filter that is briefly introduced as final part of this section.
\subsection{Discretized observation model}\label{sec:discretized_observation_model}
We start by discretizing \cref{eq:deterministic_observation_equation} for fixed $i\in\{1,\ldots,N_{obs}\}$ and fixed $j\in\{1,\ldots,N_{voxel}\}$. Numerical quadrature using a rectangular rule gives
$$c_{j,i} \approx \Delta\tau \sum_{q=0}^{N_{q}-1} c_{art}(\tau_q) \kernel_j(t_i^{obs} - \tau_q) + e_{j,i}(c_{art},\kernel)\,,$$
with $N_{q}$ equidistant abscissas $\tau_q := q \cdot \Delta\tau$ and $\Delta\tau := \frac{T}{N_{q}}$. In the original problem setting, the observation times $t_i^{obs}$ can be chosen arbitrarily. However, we here introduce a simplification, in which we assume the observation times to be given for a fixed time step size $\Delta t_{obs}$. Moreover, this time step size shall be a multiple of $\Delta \tau$, i.e. 
$$t_i^{obs} := i\, \Delta t_{obs}\,, \quad \Delta t_{obs} := s \cdot \Delta\tau\,, \quad s\in\N\,.$$
Thereby, we obtain
\begin{equation*}
\begin{aligned}
c_{j,i} &\approx \Delta\tau \sum_{q=0}^{N_{q}-1} c_{art}(\tau_q) \kernel_j(t_i^{obs} - \tau_q) + e_{j,i}(c_{art},\kernel)\\
&= \Delta\tau \sum_{q=0}^{N_{q}-1} c_{art}(q \Delta\tau) \kernel_j\left( (i\, s - q) \Delta\tau\right) + e_{j,i}(c_{art},\kernel)\,.
\end{aligned}
\end{equation*}
Since we now only need $c_{art}$ and $\kernel_j$ being evaluated at multiples of $\Delta\tau$, we can replace them by vector $\vec{c}_{art}=\left(c_{art,0},\ldots,c_{art,N_q-1}\right)^\top$ such that $c_{art,q}:=c_{art}(q \Delta \tau)$ and matrix $\mat{K}\in\R^{N_{q}\times N_{voxel}}$ with $\mat{K}:=(\kernel_{q,j})_{q,j}$ such that $\kernel_{q,j}:=\kernel_j(q \Delta\tau)$, yielding
$$c_{j,i} \approx \Delta\tau \sum_{q=0}^{N_{q}-1} c_{art,q} \, \kernel_{(i\, s -q),j} + e_{j,i}(c_{art},\kernel)\,.$$
With the extension of $k_j(t)=0$ for $t<0$ and some index substitutions, we can finally find (for each $j$, $i$) a (degenerated) matrix $\mat{H}_{j,i}\in\R^{1\times N_q}$ such that 
\begin{equation}\label{eq:deterministic_forward_map}
c_{j,i} \approx \mat{H}_{j,i} \vec{\kernel}_j + e_{j,i}(c_{art},\kernel)\,,
\end{equation}
where the $\vec{\kernel}_j\in\R^{N_q}$ are the column vectors of matrix $\mat{K}$, i.e.~$\mat{K}= (\vec{\kernel}_{0}|\ldots |\vec{\kernel}_{N_q-1})$.

Following the nomenclature of \cite{Reich2015}, we next reformulate the deterministic inference problem from \Cref{sec:full_deterministic_problem} as a sequential data assimilation problem. To this end, we first translate the involved quantities into random variables as in a Bayesian inference problem. Thereafter, we introduce the basic concepts of sequential data assimilation. 

Since the problem decouples for all voxels $j\in\{1,\ldots N_{voxel}\}$, we keep $j$ fixed for the rest of this section.

\subsection{Probabilistic view of inference}\label{sec:probabilisticViewOfInference}
Let be $(\Omega,\mathcal{F},\mathbb{P})$ a probability space. In Bayesian inference we want to gain information on a \textit{system state variable} for given observation(s). In our context, the state variable is the time-continuous kernel function $\kernel_j$. However, for simplicity and since we deal with discrete data anyway, we infer the discrete $\vec{\kernel}_j\in \R^{N_q}$ from \eqref{eq:deterministic_forward_map}, instead. Therefore, we introduce a new random variable
\begin{equation}\label{eq:system_state_random_variable}
{\vec{\mathsf{\kernel}}_j}:\Omega \rightarrow \R^{N_q}\,,
\end{equation}
replacing the time-discrete deterministic solution vector $\vec{\kernel}_j$.\footnote{We use \textsf{sans serif} letters to indicate that a given quantity is a random variable.} Moreover, we introduce a random variable
$\mathsf{e}_{j,i}(c_{art}) : \Omega \rightarrow \R,$
replacing the error term used before. Note that we assume $\vec{\mathsf{\kernel}}_j$ and $\mathsf{e}_{j,i}$ to be independent random variables. This is a rather strong simplification, since we initially modeled $e_{j,i}(c_{art},\kernel)$ to be a potentially non-linear error in the observation data, which itself is given in the indicator-dilution-theory by the arterial inflow $c_{art}$ and the tissue properties modeled by kernel $\kernel$. That is, we -- at this point -- decouple the error in the observation from the specific patient tissue. This decoupling is reflected by the new notation $\mathsf{e}_{j,i}(c_{art})$. Finally, we also consider each observation $c_{j,i}$ as random variable
$\mathsf{c}_{j,i} : \Omega \rightarrow \R,$
which is usually called \textit{observed variable}. Using \cref{eq:deterministic_forward_map}, $\mathsf{c}_{j,i}$ is defined as
\begin{equation}\label{eq:bayesian_observation_equation}
\mathsf{c}_{j,i}(\omega) := \mat{H}_{j,i} \vec{\mathsf{\kernel}}_j(\omega) + \mathsf{e}_{j,i}(c_{art})(\omega), \quad\quad \forall \omega\in\Omega\,.
\end{equation}
We will call matrix $\mat{H}_{j,i}$ \textit{(linear) forward map}.
The aim of inference is to find a reference \textit{trajectory} $\vec{\kernel}_j^{ref}$, being a realization of $\vec{\mathsf{\kernel}}_j$ such that the (measured) observations fit to the observed variable.
\subsection{Sequential data assimilation}\label{sec:sequential_data_assimilation}
Sequential data assimilation relies on an \textit{evolution model} and a \textit{forward model} to obtain $\vec{\kernel}_j^{ref}$. The models run on different time scales. The evolution model is a stochastic difference equation implying a certain predicted evolution of the system state variable over many small time steps. The forward model defines a relationship between the reference trajectory (which is to be found) and the observed data at the observation times.  
 
\subsubsection{Evolution model}\label{sec:evolutionModel}
In context of sequential data assimilation for dynamic processes, it is usually assumed that the coupling between the measured observations and the system state variable is time-local. That is, a new observation at time $t_{obs}$ only affects the system state variable for times $t\geq t_{obs}$. In our application, this is different, since the forward model, i.e.~the indicator-dilution theory, is a non-local operator in time. Therefore, we need an evolution model that allows to do global updates to the system state variable $\vec{\mathsf{k}}_j$. 
The probably most simplistic approach to model this type of global updates is given by the evolution model
\begin{equation}\label{eq:evolution_model}\vec{\mathsf{\kernel}}_j^{(l+1)} = \vec{\mathsf{\kernel}}_j^{(l)} + \sqrt{\Delta\tau}\, \vec{\mathsf{n}}^{(l)}\,,\quad l\in\{0,\ldots, N_q-1\}\,,
\end{equation}
Here $\left(\vec{\mathsf{\kernel}}_j^{(0)},\ldots,\vec{\mathsf{\kernel}}_j^{(N_q-1)}\right)$, is a sequence of random variables of the type given in \cref{eq:system_state_random_variable} for time steps $l\, \Delta\tau$. We assume $\vec{\mathsf{\kernel}}_j^{(0)} \sim \mathcal{N}(\vec{0},{\sigma_0}^2 \mat{\Sigma}_{\vec{\mathsf{n}}})$, corresponding to a zero initial guess for the kernel function with Gaussian noise with a covariance matrix $\mat{\Sigma}_{\vec{\mathsf{n}}}$. $\sigma_0\in\R$ is a scaling coefficient. The ${\vec{\mathsf{n}}}^{(l)}$s are a sequence of independent identically distributed random variables with ${\vec{\mathsf{n}}}^{(n)}:\Omega \rightarrow \R^{N_q}$ drawn as ${\vec{\mathsf{n}}}^{(l)}\sim \mathcal{N}(0,\mat{\Sigma}_{\vec{\mathsf{n}}})$.
$\mat{\Sigma}_{\vec{\mathsf{n}}}\in\R^{N_q\times N_q}$ will be chosen using a Gaussian covariance kernel such that $\mat{\Sigma}_{{\vec{\mathsf{n}}}} := \left( \sigma_{l,l^\prime} \right)_{l,l^\prime=0}^{N_q-1}$ with $\sigma_{l,l^\prime}:= \alpha e^{-\frac{\|\tau_l-\tau_{l^\prime}\|_2^2}{2 \ell^2}}$ and a parametrization in the scale $\alpha\in\R$ and the correlation length $\ell\in\R$. 

Analyzing this evolution model, we can state that it can be understood as Euler-Maruyama-based discretization of the system of stochastic ordinary differential equations
\begin{equation}\label{eq:evolution_sde}d\vec{\mathsf{\kernel}}_j = d\vec{W_t}\,,\end{equation}
where $\vec{W_t}$ is a vector of correlated univariate Wiener processes. Moreover, we observe that this evolution model, in contrast to the standard setting of dynamical processes, now only takes the role of coupling the time-discrete values in $\vec{\mathsf{\kernel}}$. This coupling is imposed by the covariance of the noise term. In fact, as we will see in Section~\ref{sec:influence_of_Xi}, the \textit{correlation length} in the Gaussian covariance kernel will have a regularizing influence on the inferred solution. Note that the choice of Gaussian noise might lead to a locally negative kernel $\mat{K}$, while this kernel is supposed to be positive. The choice of a better noise distribution is future work.

A rather natural question in context of the proposed application is, whether it would preferable to apply an EnKF-based approach for direct inversion, cf.~\cite{Iglesias2013,Schillings2017,Schillings2018} to the measurement matrix $\mat{C}$, avoiding \textit{sequential} data assimilation. In fact, we prefer sequential data assimilation, since it allows to treat the given inference problem as a time-dependent problem. In the real application case of DCE imaging, one objective of researchers is to find means to effectively control the image capturing process, in terms of a feedback loop. In that context, it is important to be able to continuously monitor the achieved approximation of the perfusion information \textit{during} the image capturing process. Based on that monitoring, one might be able to select the next observation time or the required quality for the next observation. This, however, cannot be done in direct inversion of the problem.

\subsubsection{Forward model}
We choose \cref{eq:bayesian_observation_equation} as our forward model, i.e.~we get the forward model with respect to the reference trajectory $\vec{\kernel}_j^{ref}$ 
\begin{equation}\label{eq:forward_model}
\mathsf{c}_{j,i} = \mat{H}_{j,i} {\kernel}_{j,{i\, s}}^{ref} + \mathsf{e}_{j,i}\,,\quad i\in\{1,\ldots,N_{obs}\}\,.
\end{equation}
The $\mathsf{e}_{j,i}$ are sequences of i.i.d.~random variables for growing observation time index $i$ following $\mathsf{e}_{j,i}\sim \mathcal{N}(0,{\sigma}_\mathsf{e})$, for all $i\in\{1,\ldots,N_{obs}\}$,
with $\sigma_\mathsf{e}\in\R$ the observation error variance. Note that this choice of the distribution of $\mathsf{e}_{j,i}$ is a further simplification over the simplification that has been made in Section~\ref{sec:probabilisticViewOfInference}. There, we decoupled the observation error from the kernel function $\kernel$ describing the tissue properties. Here, we further decouple the observation error from the arterial inflow and make it a purely data-independent error that is furthermore only modeled as normally distributed. It is very clear, that this is a very strong simplification. Finding a much better, maybe imaging device dependent, error is future work.

Due to $\Delta t_{obs} = s \Delta\tau$, $\kernel_{j,i\, s}^{ref}$ is the unknown reference trajectory evaluated at observation time $t_i^{obs}$.

\subsubsection{Assimilation task}
To be concise, we here only briefly summarize the general idea of the actual assimilation task with notation from \cite{Reich2015}. Further details can be found e.g.~in \cite{Reich2015}.

Let $\pi_{\vec{\mathsf{\kernel}}_j^{(i\, s)}}(\vec{\kernel}_j)$ be the probability density function of the random variable $\vec{\mathsf{\kernel}}_j^{(i\, s)}$ at time $t_i^{obs}$ for $i\in\{1, \ldots, N_{obs}\}$. Then, sequential data assimilation computes posterior PDFs
$$\pi_{\vec{\mathsf{\kernel}}_j^{(i\, s)}}(\vec{\kernel}_j|c_{j,1:i})\,,\quad i\in\{1,\ldots, N_{obs}\}\,,$$
i.e.~probability density functions of the random variables $\vec{\mathsf{\kernel}}_j^{(i\, s)}$ with an instance $\vec{\kernel}_j$ conditioned to the observations $c_{j,1}, \ldots, c_{j,i}$ that are instances of $\mathsf{c}_{j,1}, \ldots, \mathsf{c}_{j,i}$.  
This is done using an iterative approach. It is started with $\pi_{\vec{\mathsf{\kernel}}_j^{(0)}}(\vec{\kernel}_j|c_{j,1:0})$ being the PDF of $\vec{\mathsf{\kernel}}_j^{(0)}$. Then, for a given PDF $\pi_{\vec{\mathsf{\kernel}}_j^{((i-1)\, s)}}(\vec{\kernel}_j|c_{j,1:i-1})$, it iteratively 
\begin{enumerate}
\item computes the density $\pi_{\vec{\mathsf{\kernel}}_j^{\left(i\, s\right)}}(\vec{\kernel}_j|c_{j,1:i-1})$ and thereby solves a prediction problem for the given evolution model \cref{eq:evolution_model},
\item applies Bayes theorem 
$$\pi_{\vec{\mathsf{\kernel}}_j^{(i\, s)}}(\vec{\kernel}_j|c_{j,1:i})= \frac{\pi_{\mathsf{c}_{j,i}}(c_{j,i}|\vec{\kernel}_j)\, \pi_{\vec{\mathsf{\kernel}}_j^{\left(i\, s\right)}}(\vec{\kernel}_j|c_{j,1:i-1})}{\int_{\R^{N_q}}\pi_{\mathsf{c}_{j,i}}(c_{j,i}|\vec{\kernel}_j)\, \pi_{\vec{\mathsf{\kernel}}_j^{\left(i\, s\right)}}(\vec{\kernel}_j|c_{j,1:i-1})d\vec{\kernel}_j}$$
in an update step to compute $\pi_{\vec{\mathsf{\kernel}}_j^{(i\, s)}}(\vec{\kernel}_j|c_{j,1:i})$.
\end{enumerate}
 
In other words, the idea is to start from knowledge (encoded in  $\pi_{{\vec{\mathsf{\kernel}}}_j^{((i-1)\, s)}}(\vec{\kernel}_j|c_{j,1:i-1})$) at an observation time step $t_{i-1}^{obs}$. Then, knowledge for a new observation time step is forecasted / predicted using only the evolution model \cref{eq:evolution_model}. This forecast is finally corrected using the information given by observation $c_{j,i}^{}$. The unknown reference trajectory is ultimately given as mean of the marginal PDF $\pi_{\vec{\mathsf{\kernel}}_j^{(i\, s)}}(\vec{\kernel}_j|c_{j,1:N_{obs}})$.
%
%
%
%
\subsection{Ensemble Kalman Filter}
The EnKF is a Monte-Carlo--type implementation of the above discussed iterative data assimilation task. Instead of explicitly computing the posterior PDFs $\pi_{\vec{\mathsf{\kernel}}_j^{(i\, s)}}(\vec{\kernel}_j|c_{j,1:i-1})$ and $\pi_{\vec{\mathsf{\kernel}}_j^{(i\, s)}}(\vec{\kernel}_j|c_{j,1:i})$, the EnKF constructs an \textit{ensemble of realizations} of random variables representing these PDFs in an empirical sense. In that context, \textit{forecast} and \textit{analysis ensembles} are distinguished. As we will see, the computation of the forecast ensemble corresponds to approximating $\pi_{\vec{\mathsf{\kernel}}_j^{\left(i\, s\right)}}(\vec{\kernel}_j|c_{j,1:i-1,j})$, while the computation of the analysis ensemble corresponds to the approximation of $\pi_{\vec{\mathsf{\kernel}}_j^{(i\, s)}}(\vec{\kernel}_j|c_{j,1:i,j})$.

Shall $N_{e}$ be the size of the ensembles. Then, the EnKF algorithm starts by drawing $N_{e}$ samples $\vec{\kernel}_{j}^{(0),1},\ldots, \vec{\kernel}_{j}^{(0),N_{e}}$ of the (initial) system state according to the PDF of $\vec{\mathsf{\kernel}}_j^{(0)}$. The algorithm consists of two main steps which are iteratively done for $i\in\{1,\ldots,N_{obs}\}$.
\subsubsection{Forecast step}
In the forecast step, the ensemble is propagated over $s$ steps of the evolution model in \cref{eq:evolution_model} to reach the next observation time step $t_i^{obs} = i\, s\, \Delta\tau$. To achieve this, realizations $\vec{n}^{(l),m}\in\R^{N_q}$ for $m\in\{1,\ldots,N_e\}$ are drawn i.i.d.~from ${\vec{\mathsf{n}}}^{(n)}$ in each of the $s$ steps. Then the propagation equation reads for $n = 1, \ldots, s$ as 
$$ \vec{\kernel}_j^{(i(s-1)+l),m}  = \vec{\kernel}_j^{(i(s-1)+(l-1)),m} + \sqrt{\Delta\tau}\, \vec{n}^{(l),m}, \quad m\in\{1,\ldots, N_e\}\,.$$

 The newly generated ensemble is the \textit{forecast ensemble} $\left(\vec{\kernel}_{j}^{f,m}\right)_{m=1}^{N_e}$ with  $\vec{\kernel}_{j}^{f,m}:= \vec{k_j}^{i\, s,m}$. We further compute the empirical forecast mean
\begin{equation}\label{eq:empirical_forecast_mean}\overline{\vec{{\kernel}}_j^f} := \frac{1}{N_e} \sum_{m=1}^{N_e} \vec{\kernel}_{j}^{f,m} \in\R^{N_q}\end{equation}
and the empirical forecast covariance (matrix)
\begin{equation}\label{eq:empirical_forecast_covariance}\mat{\Sigma}_{\vec{\kernel}_j}^f:= \frac{1}{N_e-1} \sum_{m=1}^{N_e} \left(\vec{\kernel}_j^{f,m}-\overline{\vec{{\kernel}}_j^f}\right)\left(\vec{\kernel}_j^{f,m}-\overline{\vec{{\kernel}}_j^f}\right)^\top \in \R^{N_q\times N_q}\,.\end{equation}
\subsubsection{Analysis step}
In the analysis step, the Kalman filter \cite{Kalman1960,Reich2015} is applied to the forecast ensemble to compute an analysis ensemble $\left(\vec{\kernel}_{j}^{a,m}\right)_{m=1}^{N_e}$ representing the PDF $\pi_{\vec{\mathsf{\kernel}}_j^{(i\, s)}}(\vec{\kernel}_j|c_{j,1:i})$, which is conditioned to the new observation $c_{j,i}$. As part of the Kalman filter, the forward model \cref{eq:forward_model} with $\kernel_{j,{i\, s}}^{ref}$ being replaced by $\vec{\mathsf{\kernel}}_j^{(i\, s)}$ is evaluated. Here, we use a linear forward map $\mat{H}_{j,i}$. Moreover all involved random variables are Gaussian. Therefore, it can be shown that the analysis ensemble follows a Gaussian distribution, too and thus it can be fully characterized by the empirical analysis mean $\overline{\vec{{\kernel}}_j^a}$ and the empirical analysis covariance $\mat{\Sigma}_{\vec{\kernel}_j}^a$.

Based on this observation, the core idea of the Kalman filter is to compute the empirical analysis mean as minimization problem
$$\overline{\vec{{\kernel}}_j^a} = \argmin_{\vec{\kernel}_j\in\R^{N_q}} \frac{1}{2} \left( \left\|\vec{\kernel}_j-\overline{\vec{{\kernel}}_j^f}\right\|_{\left(\mat{\Sigma}_{\vec{\kernel}_j}^f\right)^{-1}}^2 + \left\| \mat{H}_{j,i} \vec{\kernel}_j - c_{j,i} \right\|_{{\sigma}_\mathsf{e}^{-1}}^2  \right)\,.$$
Given the linearity of $\vec{H}_{j,i}$, the minimum can be exactly computed as
$$\overline{\vec{{\kernel}}_j^a} = \overline{\vec{{\kernel}}_j^f} - \mat{U}_{j,i} (\mat{H}_{j,i} \vec{\kernel}_j - c_{j,i})\,,$$
where $U_{i,j}$ is the Kalman (update) matrix
$$\mat{U}_{j,i} = \mat{\Sigma}_{\vec{\kernel}_j^f} \mat{H}_{j,i}^\top(\mat{H}_{j,i} \mat{\Sigma}_{\vec{\kernel}_j}^f \vec{H}_{j,i}^\top + {\sigma}_{\mathsf{e}})^{-1}\,.$$
Instead of explicitly computing the empirical analysis mean and covariance (the latter by an analogous update idea), the analysis part of the Ensemble Kalman Filter (with perturbed observations) \cite[Chapter~7]{Reich2015} directly updates the forecast ensemble by
$$\vec{\kernel}_j^{a,m} = \vec{\kernel}_j^{f,m} - \mat{U}_{j,i} (\mat{H}_{j,i}\vec{\kernel}_j^{f,m}+e_{j,i,m}-c_{j,i})\,,\quad m\in\{1,\dots,N_e\}\,,$$
where $\left\{e_{j,i,m}\right\}_{m=1}^{N_e}$ are realizations of $\mathsf{e}_{j,i}$. If required, empirical versions of the analysis mean and analysis covariance can be computed analogously to \cref{eq:empirical_forecast_mean} and \cref{eq:empirical_forecast_covariance}. Finally, the next forecast step is initialized with $\vec{\kernel}_j^{(i\, s),m} = \vec{\kernel}_j^{a,m}$, that is, the analysis ensemble replaces the system state for $t_i^{obs}$.
\subsubsection{Result}
For $i=N_{obs}$ the algorithm terminates with an analysis ensemble, representing the posterior PDF\linebreak $\pi_{\vec{\mathsf{\kernel}}_j^{(N_{obs}\cdot s)}}(\vec{\kernel}_j|c_{j,1:N_{obs}})$. The reference trajectory is extracted as empirical mean ${\overline{\vec{\kernel}_j}}:= \frac{1}{N_e} \sum_{m=1}^{N_e} \vec{\kernel}_{j}^{(N_{obs}\, s),m}$. The (mean) perfusion $\overline{{p}_j}$ can be derived as $\overline{{p}_j}=\frac{1}{\rho_j}\overline{\vec{\kernel}_j}|_{t=0}$. Empirical covariances are extracted as discussed before. Moreover, in case cumulative distribution functions or other probabilistic quantities shall be extracted, a kernel-density estimator (such as \texttt{ksdensity} in Matlab) is applied to the generated ensemble. 
\section{Numerical results}\label{sec:numerical_results}
In this section, we demonstrate the beforehand introduced numerical method for artificial test data. To this end, we first introduce the source of this test data, which is a \textit{Digital Perfusion Phantom}. Then, we study the numerical properties of our method in terms of convergence, parameter dependence and input dependence in a single-voxel scenario. Finally we solve the perfusion inference problem for a slice of a full (artificial) DCE imaging brain data set.

\subsection{Digital Perfusion Phantom}
\label{sec:digital_perfusion_phantom}
\begin{figure}[th]
\centering
\scalebox{0.25}{\includegraphics{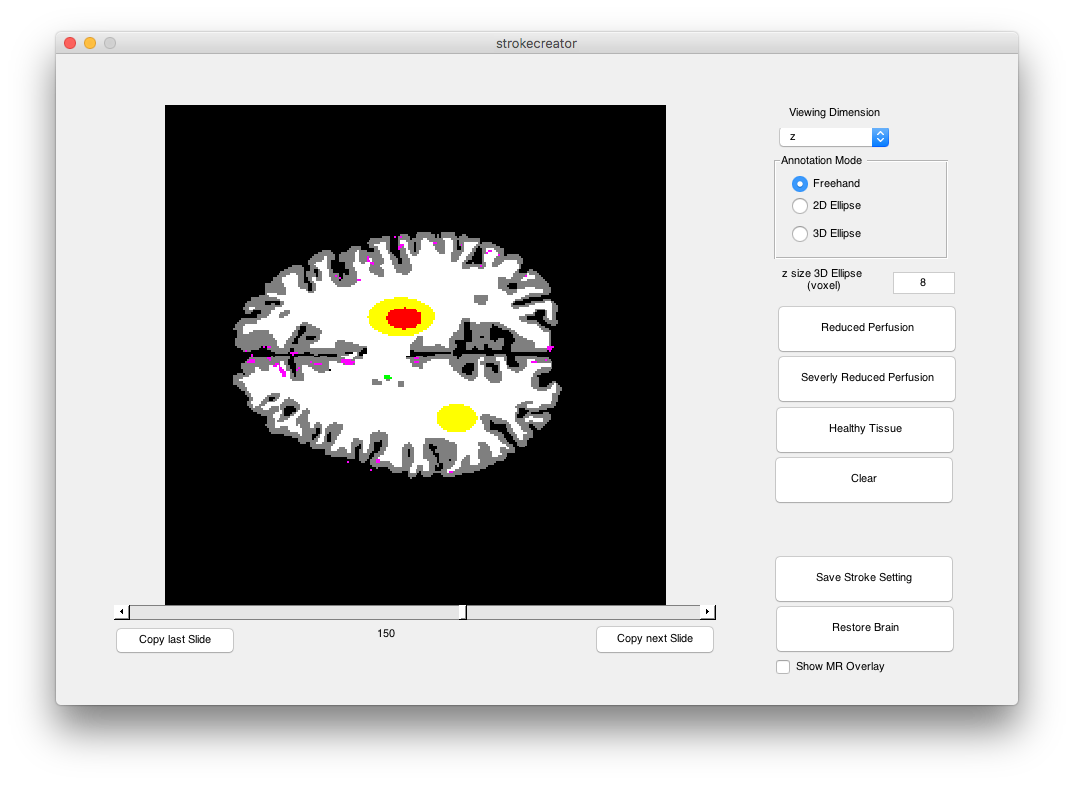}}
\caption{The source of our artificial measurements is the Digital Brain Perfusion Phantom package \cite{Manhart}. A Matlab implementation of this work is available. It allows to mark brain regions with reduced and severely reduced perfusion, here shown with the colors yellow and red. Given this data, artificial DCE imaging data is are created.}
\label{fig:screenshot_strokecreator}
\end{figure}
\begin{figure}[t]
\centering
\subfigure[arterial input function]{\label{fig:model_c_art_function}\begin{tikzpicture}[scale=0.8]
\begin{axis}[
	xlabel=time $t$,
	ylabel=$c_{art}$,
	]
\addplot+[mark=none] table {c_art.dat};
\end{axis}
\end{tikzpicture}}
\subfigure[concentration measurement in single voxel]{\label{fig:model_c_tiss_function}\begin{tikzpicture}[scale=0.8]
\begin{axis}[
	xlabel=time $t$,
	ylabel=$c_{tiss}$,
	]
\addplot+[mark=none] table {c_voi.dat};
\end{axis}
\end{tikzpicture}}
\caption{We use idealized concentration functions for one tissue voxel in order to test the implemented numerical method.}
\label{fig:model_data}
\end{figure}
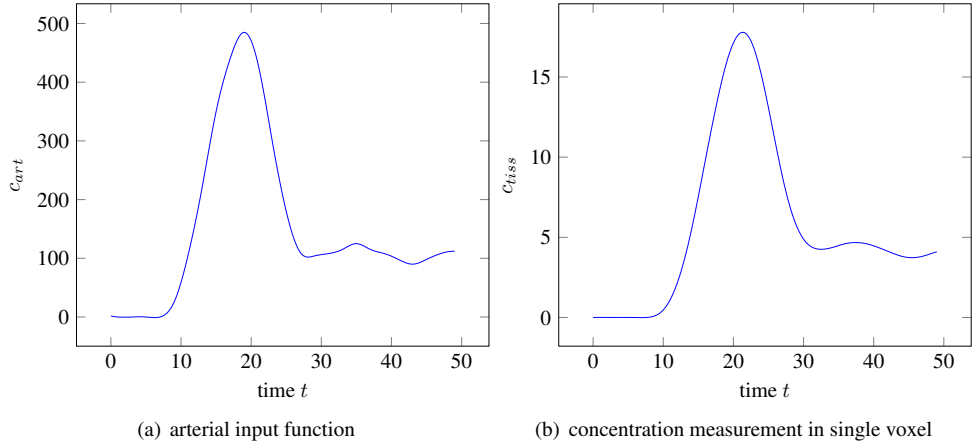

Digital Perfusion Phantoms (DPP) \cite{Riordan2011,Pianykh2012,Manhart} allow to artificially generate DCE image data for perfusion analysis. Thereby new algorithms can be tested on such data without the additional constraints of true patient data. Perfusion Phantoms basically solve the forward problem, which involves to transform perfusion information into contrast agent concentrations. In our work, we use the \textit{Digital Brain Perfusion Phantom} package \cite{Manhart}, which is a Matlab implementation of the model introduced in \cite{Riordan2011}. The software provides a radiological image of a brain. A user interface, see~\Cref{fig:screenshot_strokecreator}, allows to mark regions of reduced and strongly reduced perfusion. It is possible to control the observation snapshot time step size (i.e.~$\Delta t_{obs}$) of the artificial radiological imaging process. The measurement time is $T=49$. The resolution of the artificially generated data is $\vec{N}_\dom = (256, 256, 256)$. The arterial input function is provided as discrete evaluations $c_{art}(t_i^{input})$ with $t_i^{input} = 2\, i$. The Perfusion Phantom package uses a piecewise cubic spline interpolant through this data as exact $c_{art}$, see~\Cref{fig:model_c_art_function}. During the artificial imaging process, each snapshot (i.e.~$\vec{c}_i$) is written in a separate file. A \textit{baseline} for the radiological images is written, too. It contains the measurement data without contrast agent concentrations. In our examples, we always subtract this baseline data from the artificial measurements to obtain just the necessary concentration information.

We perform a major part of our numerical tests on a single reference voxel which has been chosen arbitrarily as $(100,130,150)$. The observation data for that single voxel is stored with a time step size of $\Delta t_{obs}=0.25$. This data is interpolated by a piecewise cubic spline to obtain measurement data at arbitrary points in time for our initial tests, cf.~\Cref{fig:model_c_tiss_function}. Towards the end of this section, results for a full slice $(\cdot,\cdot,150)$ of the full data set are discussed.  

We start by showing a series of numerical results obtained for given artificial input without noise. These results will give an insight into the choice of the different parameters of the method and into the convergence properties of the method. Noisy data is discussed afterwards.

\subsection{Data assimilation process}\label{sec:results_data_assimilation_process}
Let us first have a look at the evolution of the analysis ensemble during the sequential data assimilation process. We have chosen an observation time step size of $\Delta t_{obs}=0.25$, a quadrature step size of $\Delta\tau=0.0625$ (i.e.~$s=4$) and an ensemble size of $N_e=5000$. For a meaningful definition of the (co-)variances, we have to account for the scales of the involved quantities. By experiments, we found out that the kernel function $\kernel_j$ has a magnitude of about $10^{-3}$. Therefore, the scaling $\alpha$ of the covariance matrix $\mat{\Sigma}_{\vec{\mathsf{n}}}$ should be relative to a \textit{standard deviation} of $10^{-3}$. With this in mind, we set $\alpha = (10^{-3})^2\, 0.001$. This corresponds to a \textit{relative variance} of $0.001$. Note that it would be highly desirable to perform a coupled inference of the kernel function $\kernel_j$ \textit{and} the scaling $\alpha$. This is considered future work. The correlation length is set to $\ell=2$. For the covariance of the initial state $\vec{\mathsf{\kernel}}_j^{(0)}$, we impose an additional scaling of $\sigma_0=100$, accounting for a much larger uncertainty in the initial state. The observation error variance also needs a problem-adapted scaling. Since the concentration measurements are in the range of $10$, we shift the (co-)variance by a standard deviation of $10$. Using a relative variance of $0.0001$, we obtain ${\sigma}_\mathsf{e}=10^2\, 0.0001$.

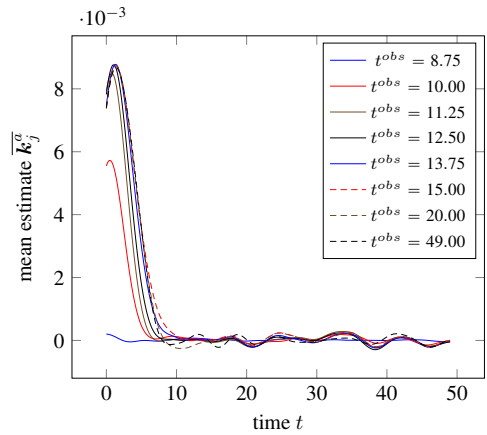
\begin{figure}[t]
\centering
\begin{tikzpicture}[scale=0.8]
\begin{axis}[
	xlabel=time $t$,
	ylabel=mean estimate $\overline{\vec{\kernel}_j^a}$,
	legend entries={$t^{obs}=8.75$,$t^{obs}=10.00$,$t^{obs}=11.25$,$t^{obs}=12.50$,$t^{obs}=13.75$,$t^{obs}=15.00$,$t^{obs}=20.00$,$t^{obs}=49.00$}
	]
\addplot+[mark=none] table {plot_r_bar_history_35.dat};
\addplot+[mark=none] table {plot_r_bar_history_40.dat};
\addplot+[mark=none] table {plot_r_bar_history_45.dat};
\addplot+[mark=none] table {plot_r_bar_history_50.dat};
\addplot+[mark=none] table {plot_r_bar_history_55.dat};
\addplot+[mark=none] table {plot_r_bar_history_60.dat};
\addplot+[mark=none] table {plot_r_bar_history_80.dat};
\addplot+[mark=none] table {plot_r_bar_history_197.dat};
\end{axis}
\end{tikzpicture}
\caption{During the sequential data assimilation process, the analysis ensemble and thereby the empirical mean of the kernel function gets continuously updated, here shown for different update time steps.}
\label{fig:data_assimilation_history_k_bar}
\end{figure}
In~\Cref{fig:data_assimilation_history_k_bar}, we show the evolution of the empirical mean $\overline{\vec{\kernel}_j^a}$, i.e.~the prediction for the unknown kernel function, for different observation times during the operation of the EnKF. Note that a scaled evolution of $\overline{\vec{\kernel}_j^a}$ at $t=0$ corresponds to the (scaled) unknown perfusion $\overline{{\vec{p}}}$. Therefore, discussing numerical results for $\overline{\vec{\kernel}_j^a}$ is equivalent to discussing results for $\overline{{\vec{p}}}$. The major information gain for the predicted result is in time interval $[10,20]$. This is the time interval in which the concentration at the arterial inlet grows. Afterwards, the data assimilation process only gains very little more information and converges towards the final result.

\subsection{Convergence in the ensemble size}
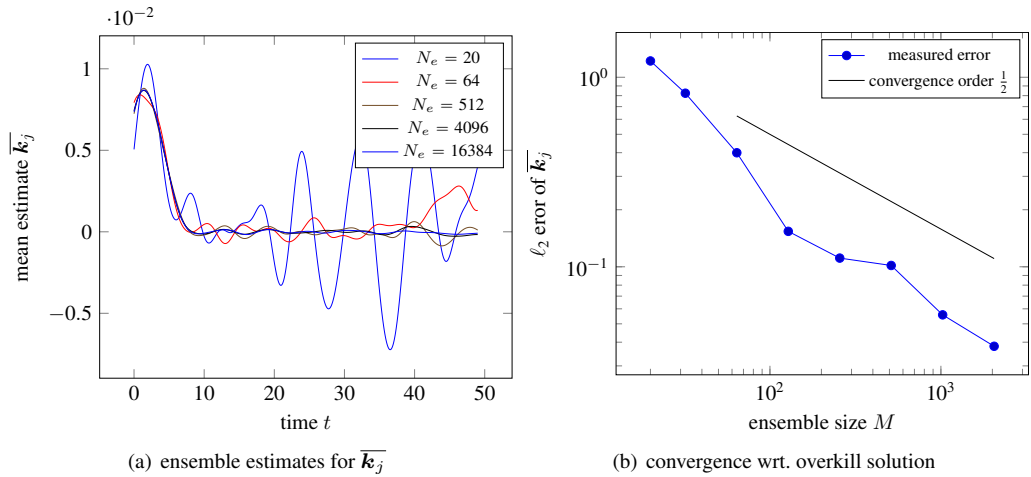
\begin{figure}[t]
\centering
\subfigure[ensemble estimates for $\overline{\vec{\kernel}_j}$]{\label{fig:convergence_ensemble_size_k_bar_plots}\begin{tikzpicture}[scale=0.8]
\begin{axis}[
	xlabel=time $t$,
	ylabel=mean estimate $\overline{\vec{\kernel}_j}$,
	legend entries={$N_e=20$,$N_e=64$,$N_e=512$,$N_e=4096$,$N_e=16384$}
	]
\addplot+[mark=none] table {plot_converence_r_bar_wrt_M_20.dat};
\addplot+[mark=none] table {plot_converence_r_bar_wrt_M_64.dat};
\addplot+[mark=none] table {plot_converence_r_bar_wrt_M_512.dat};
\addplot+[mark=none] table {plot_converence_r_bar_wrt_M_4096.dat};
\addplot+[mark=none] table {plot_converence_r_bar_wrt_M_16384.dat};
\end{axis}
\end{tikzpicture}}
\subfigure[convergence wrt.~overkill solution]{\label{fig:error_convergence_ensemble_size_k_bar}\begin{tikzpicture}[scale=0.8]
\begin{loglogaxis}[
	xlabel=ensemble size $M$,
	ylabel=$\ell_2$ eror of $\overline{\vec{\kernel}_j}$,
	legend entries={measured error, convergence order $\frac{1}{2}$}
	]
\addplot table[skip first n=0] {error_convergence_r_bar_wrt_M.dat};
\addplot [domain=64:2048] {5*x^(-1/2)};
\end{loglogaxis}
\end{tikzpicture}}
\caption{With growing ensemble size, the empirical estimate for the mean of the response / kernel function gets more accurate (left) and converges with roughly order $\frac{1}{2}$ (right).}
\label{fig:convergence_k_bar}
\end{figure}
Next, we discuss the convergence of the empirical mean of $\vec{\kernel}_j$ with respect to the ensemble size $N_e$. In the following, we will always concentrate on the last analysis ensemble obtained after assimilating the observation for $t_{obs}=49$. To shorten notation, we skip additional indices, indicating this and call the empirical mean of this analysis ensemble $\overline{\vec{\kernel}_j}$. 

Our convergence study with respect to the ensemble size uses the same parameters as in the previous paragraph. However, this time, we change the size of the ensemble. In \Cref{fig:convergence_ensemble_size_k_bar_plots}, we show the empirical mean $\overline{\vec{\kernel}_j}$ for ensemble sizes $N_e\in\{20,64,512,4096,16384\}$. The convergence in the error of the empirical mean is shown in \Cref{fig:error_convergence_ensemble_size_k_bar}. Here, we define the solution for $N_e=16384$ as overkill solution and show convergence in the relative $\ell_2$ error $\frac{\|\vec{\kernel}_j-\vec{\kernel}_j^{overkill}\|_{\ell_2}}{\|\vec{\kernel}_j^{overkill}\|_{\ell_2}}$ towards this solution. The results indicate a convergence order of approximately $\frac{1}{2}$. This is the expected order of convergence, since we use a Monte Carlo-type estimator. Note that an ensemble size of about $20$, which is often used for Ensemble Kalman Filters, seems not to be enough in this application. In that case, we observe a highly oscillatory result with a strong overshooting for the initial peak of the mean estimate (which will be the perfusion estimate).

\subsection{Convergence in the time sub-steps $\Delta\tau$}
In the following, we have a look at convergence with respect to the quadrature and evolution model step size $\Delta\tau$. Here, we do not use an overkill solution. To achieve this, we (discretely) fold the empirical mean $\overline{\vec{\kernel}_j}$ against the (discretized) arterial input function $\vec{c}_{art}$, i.e.~we transfer the prediction for $\vec{\kernel}_j$ into observation space. In observation space, we compare against the analytically given artificial measurement result $c$. Our numerical study uses a variation of the sub-step number $s$, i.e.~we change $\Delta\tau$ while keeping all other parameters as in \Cref{sec:results_data_assimilation_process}. 
\begin{figure}[t]
\centering
\subfigure[$\overline{\vec{\kernel}_j}$ in observation space]{\label{fig:plot_convergence_for_delta_tau}\begin{tikzpicture}[scale=0.8]
\begin{axis}[
	xlabel=time $t$,
	ylabel=$c_{art}\ast \overline{\vec{\kernel}_j}$,
	legend entries={$s=1$,$s=2$,$s=4$,$s=8$,$s=16$, $c_{tiss}$}
	]
\addplot+[mark=none] table {plot_converence_substepping_in_observation_space_variance_0.001_1.dat};
\addplot+[mark=none] table {plot_converence_substepping_in_observation_space_variance_0.001_2.dat};
\addplot+[mark=none] table {plot_converence_substepping_in_observation_space_variance_0.001_4.dat};
\addplot+[mark=none] table {plot_converence_substepping_in_observation_space_variance_0.001_8.dat};
\addplot+[mark=none] table {plot_converence_substepping_in_observation_space_variance_0.001_16.dat};
\addplot+[mark=none] table {plot_converence_substepping_in_observation_space_reference_variance_0.001.dat};
\end{axis}
\end{tikzpicture}}
\subfigure[convergence]{\label{fig:error_convergence_for_delta_tau}\begin{tikzpicture}[scale=0.8]
\begin{loglogaxis}[
	xlabel=time step size $\Delta\tau$,
	ylabel=error in $\vec{c}_{art} \ast \overline{\vec{\kernel}_j}$ wrt.~$c$,
	legend entries={measured error, convergence order $1$, convergence order $\frac{1}{2}$},
	legend pos=south east
	]
\addplot table[x index=0,y index=1] {error_convergence_substepping_in_observation_space_variance_0.001.dat};
\addplot+ [no marks, dashed, domain=0.02:0.2] {0.07*x^(1)};
\addplot+ [no marks, domain=0.02:0.2] {0.11*x^(1/2)};
\end{loglogaxis}
\end{tikzpicture}}
\caption{A smaller time step size for the quadrature / system state model leads to convergence of $\overline{\vec{\kernel}_j}$ in observation space towards the measurement concentration $c$.}
\label{fig:convergence_M}
\end{figure}
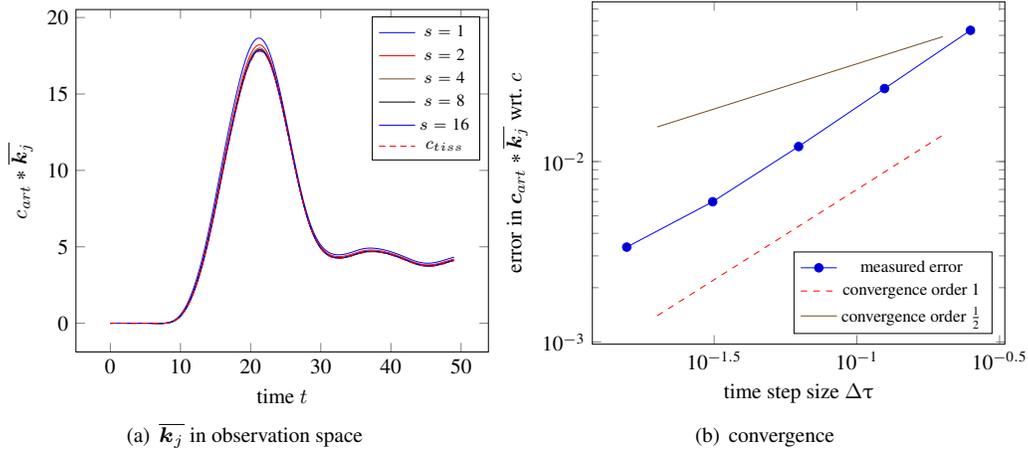

 In~\Cref{fig:plot_convergence_for_delta_tau}, we visually compare the results obtained for an increasing number of sub-steps $s$ (i.e.~decreasing $\Delta\tau$). The convergence plot in~\Cref{fig:error_convergence_for_delta_tau} further shows the error reduction in the relative $\ell_2$ norm for decreasing $\Delta\tau$ if we compare the convolved mean estimate $\overline{\vec{\kernel}_j}$ with the real observation data. The results indicate approximately first order convergence. In fact, parameter $\Delta\tau$ influences the Euler-Maruyama approximation of the continuous stochastic differential equation \cref{eq:evolution_sde} and the quadrature of the convolution integral. While the Euler-Maruyama method is known to have halve order convergence, the rectangular rule is convergent of second order for sufficiently smooth integrands. The observed convergence behavior strongly depends on the dominance of one of the errors (time-integration, quadrature). The observed first order seems to indicate that the quadrature error for the convolution integral is dominant. Nevertheless, full second order convergence is not achieved. This observation is clearly a pre-asymptotic and strongly problem-dependent result.

\subsection{Influence of the correlation length in the system state noise}\label{sec:influence_of_Xi}
\begin{figure}[t]
\centering
\begin{tikzpicture}[scale=0.8]
\begin{axis}[
	xlabel=time $t$,
	ylabel=mean estimate $\overline{\vec{\kernel}_j}$,
	legend entries={$\ell=0.125$,$\ell=0.5$, $\ell=2$}
	]
\addplot+[mark=none] table {plot_r_bar_wrt_correltion_lengths_0.125000.dat};
\addplot+[mark=none] table {plot_r_bar_wrt_correltion_lengths_0.500000.dat};
\addplot+[mark=none] table {plot_r_bar_wrt_correltion_lengths_2.000000.dat};
\end{axis}
\end{tikzpicture}
\caption{Longer correlation lengths $\ell$ impose a higher smoothness on the ensemble estimate.}
\label{fig:influence_of_correlation_length}
\end{figure}
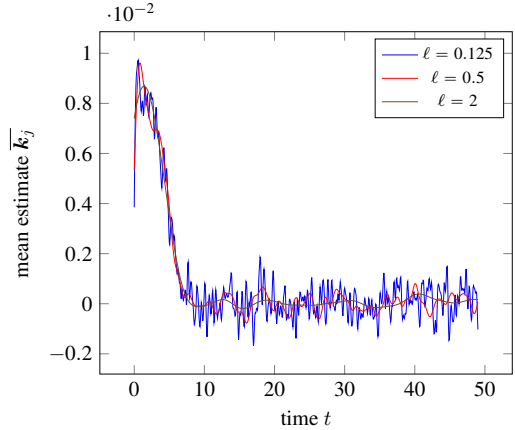

Our next study shall give an insight into the influence of the system state model, more specifically the influence of the correlation length $\ell$ of the random variable ${\vec{\mathsf{n}}}$ on the inferred solution. To study the influence of the correlation length, we keep the parameters as in \Cref{sec:results_data_assimilation_process} and apply different correlation lengths $\ell\in\{0.125, 0.5, 2\}$. The results of this numerical study are given in \Cref{fig:influence_of_correlation_length}. Here, the inferred kernel function $\overline{\vec{\kernel}_j}$ is shown for different correlation lengths. For growing correlation length the result gets less noisy. Hence, a larger correlation length has a regularizing effect on the solution. Since, in general, we seek for smooth solutions, we always choose $\ell=2$.

\subsection{Influence of the number of observations}
\begin{figure}[t]
\centering
\begin{tikzpicture}[scale=0.8]
\begin{axis}[
	width=0.6\textwidth,
	height=0.475\textwidth,
	xlabel=$\kernel|_{t=0}$,
	ylabel=estimated PDF for $\kernel|_{t=0}$,
	legend entries={$\Delta t_{obs}=1.0$,$\Delta t_{obs}=0.5$,$\Delta t_{obs}=0.25$,$\Delta t_{obs}=0.125$},
	legend pos=north west
	]
\addplot+[mark=none] table {plot_pdf_wrt_delta_t_obs_1.000000e+00.dat};
\addplot+[mark=none] table {plot_pdf_wrt_delta_t_obs_5.000000e-01.dat};
\addplot+[mark=none] table {plot_pdf_wrt_delta_t_obs_2.500000e-01.dat};
\addplot+[mark=none] table {plot_pdf_wrt_delta_t_obs_1.250000e-01.dat};
\end{axis}
\end{tikzpicture}
\caption{The more observation samples are taken, the more reliable the estimate of the solution. Hence, the estimated PDF for $\kernel|_{t=0}$ shows a smaller variance for smaller observation time steps $\Delta_{t_{obs}}$.}
\label{fig:influcence_of_number_of_observations}
\end{figure}
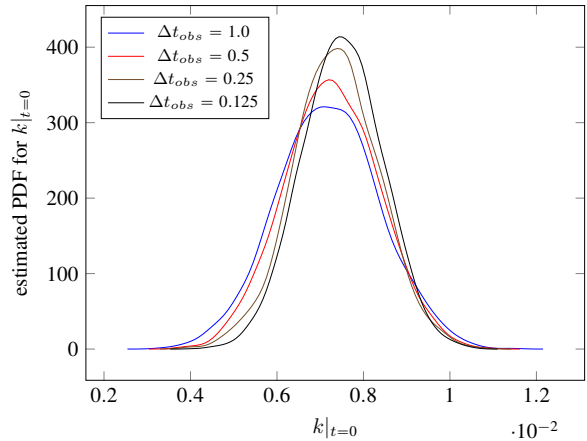
Our final test with noise-free model data on a single voxel highlights the influence of a change of the observation time step size $\Delta t_{obs}$, i.e.~a change in the number of observations that are made during the imaging process. To test this, we take the same parameters as in \Cref{sec:results_data_assimilation_process}, but change the observation time step size as $\Delta t_{obs}\in\{0.125, 0.25, 0.5, 1.0\}$ while keeping $\Delta\tau$ constant. The quantity that we study is the computed probability density function for $\kernel|_{t=0}$, hence a scaled version of $\overline{{\vec{p}}}$. We use the kernel density estimator \texttt{ksdensity} in Matlab to reconstruct a continuous PDF for the ensemble data.   

The results of this study can be seen in \Cref{fig:influcence_of_number_of_observations}. Here, we make two observations. First, the mean of the PDF still changes for growing number of measurements, converging towards a true solution. Second, and more important, we observe a variance reduction if we increase the number of measurements. This type of information would not be available in classical inverse approaches for compute perfusion estimation. That is, we can now obtain confidence information for our solution.

\subsection{Inference from noisy data}
Until now, we considered noise-free input data. Instead, we now discuss the same one-voxel input as before, but add artificial noise as
$${c}_{j,i}^{noisy} = c_{j,i} + w_{j,i}\,,\quad i\in\{1,\ldots, N_{obs}\}\,,$$
where the $w_{j,i}$ are realizations of i.i.d.~random variables $\mathsf{w}_{j,i}:\Omega\rightarrow\R,\,\mathsf{w}_{j,i}\sim\mathcal{N}(0,{\sigma}_{\mathsf{w}})$ with ${\sigma}_{\mathsf{w}}\in\R$ the variance of the noise.

\begin{figure}[t]
\centering
\subfigure[]{\label{fig:plot_results_with_noisy_data_0.25}\begin{tikzpicture}[scale=0.8]
\begin{axis}[
	xlabel=time $t$,
	ylabel=contrast agent concentration,
	ymax=30,
	ymin=-2,
	legend entries={$\mathbf{c}_j^{noisy}$ for $\alpha_{rel}^2\approx 0.0156$,$c_{art}\ast \overline{\vec{\kernel}_j}$},
	legend pos=north east
	]
\addplot+[mark=none] table[x index=0,y index=1] {plot_noisy_results_wrt_variance_1.562500.dat};
\addplot+[mark=none] table[x index=0,y index=2] {plot_noisy_results_wrt_variance_1.562500.dat};
\end{axis}
\end{tikzpicture}}
\subfigure[]{\label{fig:plot_results_with_noisy_data_1.0}\begin{tikzpicture}[scale=0.8]
\begin{axis}[
	xlabel=time $t$,
	ylabel=contrast agent concentration,
	ymax=30,
	ymin=-2,
	legend entries={$\mathbf{c}_j^{noisy}$ for $\alpha_{rel}^2=0.0625$,$c_{art}\ast \overline{\vec{\kernel}_j}$},
	legend pos=north east
	]
\addplot+[mark=none] table[x index=0,y index=1] {plot_noisy_results_wrt_variance_6.250000.dat};
\addplot+[mark=none] table[x index=0,y index=2] {plot_noisy_results_wrt_variance_6.250000.dat};
\end{axis}
\end{tikzpicture}}
\caption{Even for stronger noise on the input data the inference of $\overline{\vec{\kernel}_j}$ is acceptable, as long as the measurement variance ${\sigma}_{\mathsf{e}}$ is chosen appropriately. Here, we compare the noisy input $\vec{c_j}^{noisy}$ $\overline{\vec{\kernel}_j}$ in observation space for $\alpha_{rel}^2=0.015625$ (left) and $\alpha_{rel}^2=0.0625$ (right).}
\label{fig:results_with_noisy_data}
\end{figure}
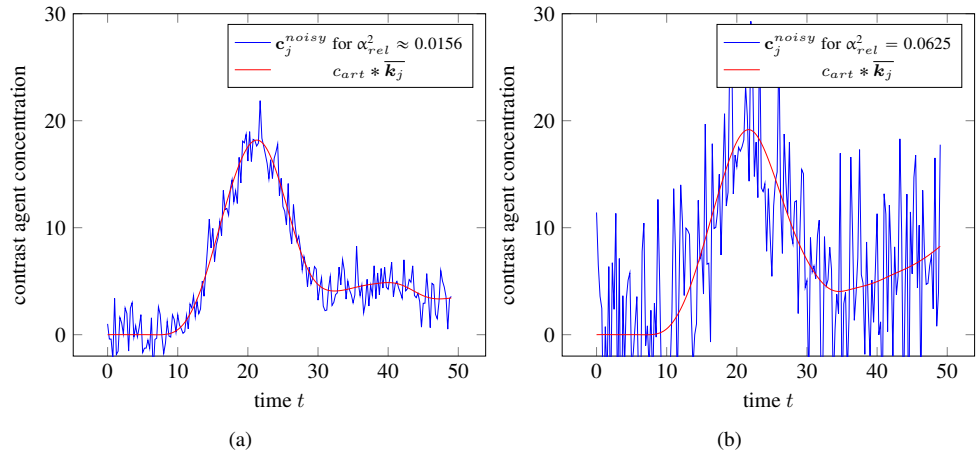

We use a series of test cases with ${\sigma}_{\mathsf{w}}=10^2\alpha_{rel}$ and $\alpha_{rel}\in\{2^{-10},2^{-8},2^{-6},2^{-4},2^{-2}\}$. Hence, $\alpha_{rel}$ corresponds to the relative variance with respect to the magnitude of the measurements. We keep a major part of the parameters from \Cref{sec:results_data_assimilation_process}. However, we change the fixed observation error variance ${\sigma}_{\mathsf{e}}$, to a problem-adapted one, namely, i.e.~${\sigma}_{\mathsf{e}}={\sigma}_{\mathsf{w}}$. Note that in practice, one would empirically \textit{estimate} the noise in the measurement data and would set ${\sigma}_{\mathsf{e}}$ accordingly. Another change concerns the number of samples $N_e$ in the EnKF. As our experiments showed, the size of the ensemble has to be increased for higher variances. This is well covered by classical Monte-Carlo theory. Therefore, we set $N_e=10000$ for $\alpha_{rel}\in\{2^{-10},2^{-8},2^{-6}\}$ while we use $N_e = 60000$ and $N_e=100000$ for $\alpha_{rel}=2^{-4}$ and $\alpha_{rel}=2^{-2}$, respectively.

\begin{figure}[t]
\centering
\begin{tikzpicture}[scale=0.8]
\begin{axis}[
	xlabel=$\vec{\bar{\kernel}}_j$,
	ylabel=estimated PDF for $\kernel|_{t=0}$,
	legend entries={$\alpha_{rel}^2\approx 0.0001$,$\alpha_{rel}^2\approx 0.0039$,$\alpha_{rel}^2\approx 0.0156$,$\alpha_{rel}^2=0.0625$,$\alpha_{rel}^2=0.25$},
	xmin=-0.003,
	xmax=0.015,
	legend pos=north west
	]
\addplot+[mark=none] table[x index=0,y index=1] {plot_pdf_of_noisy_results_wrt_variance_0.097656.dat};
\addplot+[mark=none] table[x index=0,y index=1] {plot_pdf_of_noisy_results_wrt_variance_0.390625.dat};
\addplot+[mark=none] table[x index=0,y index=1] {plot_pdf_of_noisy_results_wrt_variance_1.562500.dat};
\addplot+[mark=none] table[x index=0,y index=1] {plot_pdf_of_noisy_results_wrt_variance_6.250000.dat};
\addplot+[mark=none,color=green] table[x index=0,y index=1] {plot_pdf_of_noisy_results_wrt_variance_25.000000.dat};
\end{axis}
\end{tikzpicture}
\caption{The proposed method is pretty robust with respect to noise. This can be seen, if we study the estimated probability density functions for $\kernel|_{t=0}$. With growing noise variance, the empirical PDF estimate still recovers the mean appropriately. Extreme noise variances degenerate the result, as expected.}
\label{fig:pdf_comparison_for_noisy_data}
\end{figure}
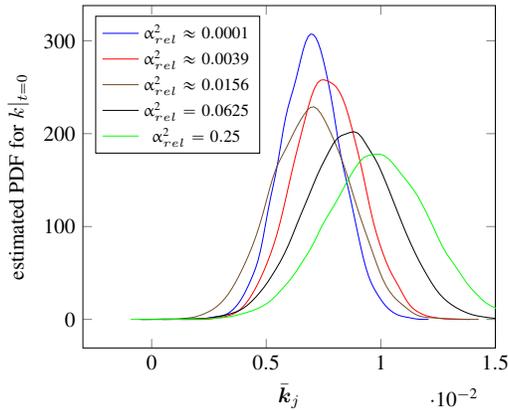

In \Cref{fig:results_with_noisy_data}, we give two examples of noisy inputs for $\alpha_{rel}^2=0.015625$ and $\alpha_{rel}^2=0.0625$. In the latter case, the original input signal is already severely degenerated. The predicted mean solutions in observation space $\vec{c}_{art}\ast \overline{\vec{\kernel}_j}$ are also given in~\Cref{fig:results_with_noisy_data}. In fact, the reconstructed solution is almost not influenced for $\alpha_{rel}^2=0.015625$ and gets a little distorted for stronger noise. In \Cref{fig:pdf_comparison_for_noisy_data}, we compare estimates of the PDF for $\kernel|_{t=0}$ for growing noise in the data. Since we appropriately account for the noise in the input, the mean is almost identical up to $\alpha_{rel}^2=0.015625$. For higher relative noise variances, the probability density functions still cover the general tendency of the results. Note that the variance in the solutions grows for larger noise in the input. This effect is not primarily caused by the noisy input, but by the imposed observation error $\mathsf{e}_{j,i}$, which acts here as a regularization for the noisy input. Nonetheless, as long as the observation error variance is set in the range of the input noise variance, the variance in the solution correctly represents the variance coming from the noise in the input.

\subsection{Application problem}
We finally apply the beforehand studied method to a full application problem given by the Digital Brain Perfusion Phantom introduced in \Cref{sec:digital_perfusion_phantom}. We use the slice $(\cdot, \cdot, 150)$ with the choice of regions with reduced and severely reduced perfusion as in~\Cref{fig:screenshot_strokecreator}. We discuss a result for a large observation time-step size $\Delta t_{obs}=1.0$, a highly resolved quadrature with $\Delta\tau = 0.0625$ and a noise with variance ${\sigma}_{\mathsf{w}}=10^2\, 0.015625$. The observation error variance is adapted as ${\sigma}_{\mathsf{e}}={\sigma}_{\mathsf{w}}$. All other parameters are kept as in \Cref{sec:results_data_assimilation_process}.

\subsubsection{Storage and performance considerations}
Storing and computing the ensembles for the discussed test cases is a rather challenging task. Just considering the analysis ensemble for a single slice, we need to store for each of the $256\times 256$ voxels $5000$ realizations of discrete kernel functions $\vec{\kernel}_j$ given via $N_q=785$ double precision values leading to a total storage requirement of
$$ 256\times 256 \times 5000 \times 785 \times 8\ \mbox{Bytes} \approx 1834\ \mbox{GBytes}\,.$$

All our calculations are done in Matlab. We always compute $8$ rows of the final $256\times 256$ slice at the same time and reuse the random input for each voxel in order to reduce the runtime. Note that especially sampling from ${\vec{\mathsf{n}}}^{(n)}$ is very computationally demanding. In order to do the calculations, we need constant access to way more than 64 GBytes of RAM. Due to storage und memory requirements, we use nodes of the cluster \textit{Rhea} at Oak Ridge National Lab to compute the full problem. Each node has 128 GBytes of RAM and a dual Intel® Xeon® E5-2650 CPU with 16 cores. To compute $8$ lines, i.e.~results for $8\times 256=2048$ voxels, we need about $3$ hours and $15$ minutes, noting that Matlab uses approximately 14 cores of the full machine.  The total computing time (with respect to one node of Rhea) is thereby roughly $104$ hours or about $4.3$ days on a single machine.

Even though this amount of computing time seems to be rather prohibitive for the specific application case, it is clear that the discussed algorithm is extremely easy to parallelize. Especially, it seems to be very well suited to a parallelization on graphics processing units (GPUs) or other many-core hardware, as long as the results of the calculation are constantly streamed out to CPU memory. An appropriate parallel implementation is future work.

\pgfplotsset{
        colormap={parula}{
            rgb255=(53,42,135)
            rgb255=(15,92,221)
            rgb255=(18,125,216)
            rgb255=(7,156,207)
            rgb255=(21,177,180)
            rgb255=(89,189,140)
            rgb255=(165,190,107)
            rgb255=(225,185,82)
            rgb255=(252,206,46)
            rgb255=(249,251,14)
        },
	colormap={low_perfusion}{rgb255(0cm)=(255,255,255); rgb255(1cm)=(53,42,135)},
	colormap={medium_perfusion}{rgb255(0cm)=(255,255,255); rgb255(1cm)=(21,177,180)},
	colormap={high_perfusion}{rgb255(0cm)=(255,255,255); rgb255(1cm)=(225,185,82)}
}
\begin{figure}[th]
\centering
\subfigure[Approximated perfusion ${\overline{\vec{p}}}$\label{fig:full_approximated_perfusion}]{
\begin{tikzpicture}[scale=0.75]
\begin{axis}[
view={0}{90},
xlabel=$x$,
ylabel=$y$,
colormap name=parula,
colorbar
]
\addplot3[surf,shader=interp,mesh/rows=256] table[col sep=space]{perf_tilde_real_all.dat};
\end{axis}
\end{tikzpicture}}
\subfigure[Noise-free reference perfusion given by the DPP\label{fig:full_perfusion_reference_solution}]{\begin{tikzpicture}[scale=0.75]
\begin{axis}[
view={0}{90},
xlabel=$x$,
ylabel=$y$,
colormap name=parula,
colorbar,
zmax=1.0
]
\addplot3[surf,shader=interp,mesh/rows=256] table[col sep=space]{perf_reference_all.dat};

\end{axis}
\end{tikzpicture}}

\caption{Our approximation method recovers the reference perfusion result (right) as mean of the ensemble in the Ensemble Kalman Filter. Both results match well, even though we introduced a considerable amount of artificial noise.}

\end{figure}
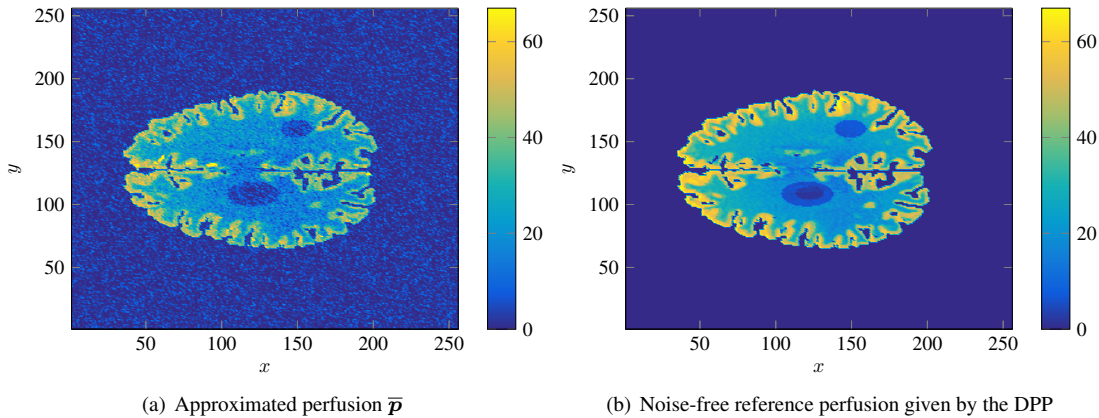

\subsubsection{Quantities of interest}\label{sec:quantities_of_interest}
{In our application examples, we consider the approximation of probabilistic quantities of interest in connection with the perfusion $p_j:=p(\kernel_j) = \frac{1}{\rho_j} \kernel(0)$. Besides of the mean $\overline{p_j}$ we are especially interested in probabilities for the corresponding random variable ${\sf p}_j$ to be in a given range. To be more specific, we compute the probabilities
$$\mathbb{P}({\sf p}_j < 10)\,,\quad \mathbb{P}(20\leq {\sf p}_j < 40)\,,\quad \mathbb{P}({\sf p}_j \geq 50)\,,$$
noting that the underlying perfusion $p_j$ lies in the interval $[0,70]$ in the case of the Digital Brain Perfusion Phantom data that we consider. These quantities give probabilities for low, medium and high perfusion in some region of the brain. Given the final analysis ensemble for $\kernel_j$, it is easy to compute the above quantities by using the kernel density estimator \texttt{ksdensity}. The latter one can compute a cumulative distribution function (CDF) for each voxel, which is finally evaluated appropriately.

\subsubsection{Discussion of results}
\begin{figure}[th]
\centering
\begin{tikzpicture}[scale=0.9]
\begin{axis}[
view={0}{90},
xlabel=$x$,
ylabel=$y$,
colormap name=low_perfusion,
colorbar
]
\addplot3[surf,shader=interp,mesh/rows=256] table[col sep=space]{prob_cbf_smaller_10_real_tilde_all.dat};

\end{axis}
\end{tikzpicture}\vspace*{-0.5em}
\caption{The advantage of the proposed method is that we are now also able to compute probabilistic information for the solution, here shown by plotting the probability $\mathbb{P}({\sf p}_j < 10)$. Hence, the depicted results give the space-dependent probability for low ($<10$) perfusion.}
\label{fig:probability_low_perfusion}
\end{figure}
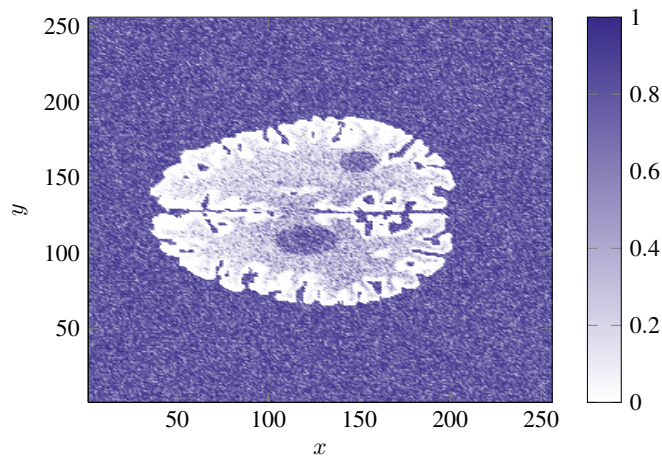

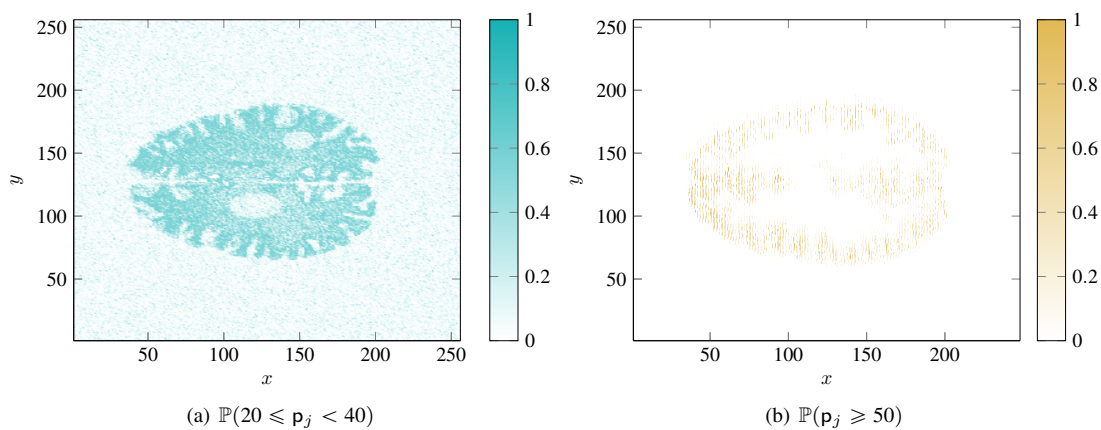
\begin{figure}[th]
\centering
\subfigure[$\mathbb{P}(20\leq {\sf p}_j < 40)$\label{fig:probability_medium_perfusion}]{\begin{tikzpicture}[scale=0.75]
\begin{axis}[
view={0}{90},
xlabel=$x$,
ylabel=$y$,
zmax=1.0,
colormap name=medium_perfusion ,
colorbar
]
\addplot3[surf,shader=interp,mesh/rows=256] table[col sep=space]{prob_cbf_between_20_40_real_tilde_all.dat};

\end{axis}
\end{tikzpicture}}
\subfigure[$\mathbb{P}({\sf p}_j \geq 50)$\label{fig:probability_high_perfusion}]{\begin{tikzpicture}[scale=0.75]
\begin{axis}[
view={0}{90},
xlabel=$x$,
ylabel=$y$,
zmax=1.0,
colormap name=high_perfusion,
colorbar
]
\addplot3[surf,shader=interp,no marks,mesh/rows=256] table[col sep=space]{prob_cbf_larger_50_real_tilde_all.dat};

\end{axis}
\end{tikzpicture}}
\caption{Based on the results of the EnKF, it is easily possible to identify regions of high probability to have medium (left) and high (right) perfusion.}
\label{fig:}
\end{figure}

An important advantage of the use of a Digital Perfusion Phantom is the existence of a reference solution to compare with. The DPP software that we use stores the reference solution together with the other generated data. In \Cref{fig:full_perfusion_reference_solution}, we show the reference solution for our full application test case. The approximated result of our application example study, i.e. {$\overline{\vec{p}}$}, is shown in \Cref{fig:full_approximated_perfusion}. As expected from our single-voxel study, the inferred perfusion matches the exact perfusion result well. Note that this is the case even though we add a considerable amount of noise on the measurements. 

As discussed before, we can use the ensemble-based estimate of the posterior probability density function to extract a wide range of probabilistic information on the inferred solution. This is the main result of this work. To exemplify this, we compute space-dependent probabilities for low (\Cref{fig:probability_low_perfusion}), medium (\Cref{fig:probability_medium_perfusion}) and high (\Cref{fig:probability_high_perfusion}) perfusion ranges, cf.~\Cref{sec:quantities_of_interest}. In case of \Cref{fig:probability_low_perfusion}, we e.g.~can now easily identify ranges of low perfusion and even can give a probability for this result.

In general, we claim that this probability information or derived probabilistic quantities (variance, percentiles, etc.) can give domain-experts in radiology a much clearer information on the reliability of the inferred estimates.

\section{Summary}\label{sec:summary}
In this work, we have discussed the use of Ensemble Kalman Filters for sequential data assimilation in order to infer probabilistic information on (blood) perfusion in tissue for given measurements from dynamic contrast--enhanced imaging. The deterministic inference of perfusion is well-known in the field of radiological imaging. However, to the best of the author's knowledge, the new contribution is the approximation of PDFs for the perfusion given (noisy) measurements. EnKF are well-known in inference for dynamical systems and partial differential equations with stochastic coefficients. Hence, modeling the dynamic contrast--enhanced imaging process as sequential data assimilation in a Bayesian context was the main contribution of the work. Given the ensemble-based approximation of the PDF, we could compute probabilistic quantities such as probabilities for perfusion parameter ranges. 

The new approach was first investigated for a single-voxel example with respect to convergence and parameter influence. Afterwards, it was applied to artificial application data generated by a Digital Perfusion Phantom, i.e.~a model for deriving DCE image data for given perfusion data. Overall, the effectiveness of the method could be demonstrated, showing 
empirical convergence results and appropriate approximations of probabilistic information. The use of realistic patient data, refined problem-adapted covariance kernels, advanced filtering techniques and an efficient parallel implementation are future work.

\section*{Acknowledgements}
This work is funded by the Swiss National Science Foundation (SNF) under project number $407540\_167186$. Furthermore, this research used resources of the Oak Ridge Leadership Computing Facility at the Oak Ridge National Laboratory, which is supported by the Office of Science of the U.S.~Department of Energy under Contract No.~DE-AC05-00OR22725.

The author also likes to thank Wolfram Stiller and Christian Weis of the department of \textit{Diagnostic and Interventional Radiology} of the \textit{University Medical Center Heidelberg} and Holger Fr\"oning of the \textit{Institute of Computer Engineering at University of Heidelberg} for fruitful initial discussions on the application background.

\bibliographystyle{IJ4UQ_Bibliography_Style}
\bibliography{perfusion_uq}

\end{document}